\documentclass[12pt]{article}%28.08.2017
\usepackage{amsmath}
\usepackage{amssymb}
\usepackage{dsfont}
\usepackage{cite}
\newsymbol \blackbox 1004
\newcommand{\eh}{\hfill}\newlength{\sperr}

\newenvironment{proof}{{\settowidth{\sperr}{\bf\rm
Proof}%
\par\addvspace{0.3cm}\noindent\parbox[t]{1.3\sperr}
{\bf\rm P\eh r\eh o\eh o\eh f\eh }%
}}{\nopagebreak\mbox{}
$\blackbox$\par\addvspace{0.3cm}}

\def\nn{\nonumber}
\def\a{\alpha}
\def\b{\beta}
\def\g{\gamma}

\def\vk{\varkappa}

\def\s{\sigma}
\def\la{\lambda}
\def\om{\omega}

\def\vp{\varphi}
\def\vt{\vartheta}

\def\wh{\widehat}
\def\wt{\widetilde}
\def\ov{\overline}

\def\BC{{\mathbb C}}
\def\BR{{\mathbb R}}
\def\BN{{\mathbb N}}
\def\clp{{\mathcal P}}

\def\cld{{\mathcal D}}
\def\cle{{\mathcal E}}

\def\clh{{\mathcal H}}

\def\clk{{\mathcal K}}

\def\clq{{\mathcal Q}}

\def\clw{{\mathcal W}}
\def\cld{{\mathcal D}}

\def\mfa{{\mathfrak A}}

\newcommand{\E}{\mathrm{e}}
\newcommand{\I}{\mathrm{i}}

\newtheorem{Pa}{Paper}[section]
\newtheorem{Tm}[Pa]{{\bf Theorem}}

\newtheorem{Cy}[Pa]{{\bf Corollary}}
\newtheorem{Rk}[Pa]{{\bf Remark}}

\newtheorem{Dn}[Pa]{{\bf Definition}}

\newtheorem{Pn}[Pa]{{\bf Proposition}}

\title{Verblunsky-type coefficients \\ for Dirac and canonical systems  \\ generated by Toeplitz and Hankel matrices,
respectively}

\author{Alexander Sakhnovich}

\date{}
\begin{document}
\maketitle

\begin{abstract} We introduce Verblunsky-type coefficients of Toeplitz and Hankel matrices, which correspond
to the discrete Dirac and canonical systems generated by Toeplitz and Hankel matrices,
respectively. We prove one to one correspondences between positive-definite Toeplitz (Hankel) matrices
and their Verblunsky-type coefficients as analogs of the well-known Verblunsky's theorem.
Several interconnections with the spectral theory are described as well.

\end{abstract}

\vspace{0.3em}


\vspace{0.3em}


\section{Introduction}\label{Intro}
\setcounter{equation}{0}
Orthogonal scalar and matrix polynomials is an important and actively studied domain
(see, e.g., \cite{Ak, Atk, Beals, DaS, DFK, Sz, Khru, KiSi, Krein, KrS, SiOPUC, Van} and numerous references therein).
In his interesting talk on the spectral theory of orthogonal polynomials \cite{Sitalk},
B.~Simon mentioned ``two tribes" working in this domain. Namely, the specialists
in orthogonal polynomials and the specialists in spectral theory, which communicate
quite insufficiently although their results are closely connected. We would like
to discuss here  one more set of related results, which could be fruitfully
used in the theory of orthogonal polynomials and corresponding problems of Toeplitz 
and Hankel matrices.

We recall that the theory of orthogonal polynomials on the unit circle (OPUC)
studies interrelations between a measure $d \tau$ (or, equivalently, nondecreasing weight function
$\tau(t)$ on $[-\pi, \, \pi]$) and positive-definite Toeplitz matrices
\begin{align}\label{01}
S(n)=\{s_{j-i}\}_{i,j=1}^n, \quad s_k=\frac{1}{2 \pi}\int_{-\pi}^{\pi}\E^{\I k t}d\tau(t)
\end{align}
on one side,  and orthogonal trigonometric polynomials $\clp_r(\E^{\I t})$ generated by the measure $d \tau$ on the other
side. The theory of orthogonal polynomials on the real axis (OPRL) studies interrelations between a  nondecreasing weight function
$\tau(t)$ on $(-\infty, \, \infty)$ and Hankel matrices
\begin{align}& \label{02}
H(n)=\{H_{i+j-2}\}_{i,j=1}^n, \quad H_k=\int_{-\infty}^{\infty}t^{k}d\tau(t)
\end{align}
on one side,  and orthogonal polynomials generated by the measure $d \tau$ on the other
side. 

It is well known (see, e.g., \cite[Sect. 5.2]{Ak}) that the orthornormal  polynomials $\clp_r(\la)$ on the unit circle
satisfy Szeg\"o recurrence
\begin{equation} \label{03}
Z_{k+1}(\la)=\frac{1}{\sqrt{1-|a_k|^2}}\begin{bmatrix} 1 &-\ov{a_k} \\ - {a_k} & 1
\end{bmatrix} \left[
\begin{array}{cc}
\la I_{p} & 0 \\ 0 & I_{p}
\end{array}
\right] Z_k(\la), 
\end{equation}
where $Z_r(\la):=\left[
\begin{array}{c} \clp_r(\la) \\ [8pt] \la^r \ov{\clp_r(1/\ov{\la})}
\end{array}
\right]$, and the coefficients $\{a_k\}$ are so called Verblunsky coefficients:
\begin{align}\label{04}
|a_k|<1 \quad (0 \leq k <\infty).
\end{align}
Clearly, there is a one to one correspondence between Szeg\"o recurrences \eqref{03}, \eqref{04} and the sequences of Verblunsky
coefficients $\{a_k\}$. A fundamental Verblunsky theorem (see, e.g., a detailed discussion in \cite{SiOPUC}) states
that there is a one to one correspondence between Verblunsky coefficients and nontrivial probability measures $d \tau$
on $[-\pi, \, \pi]$. 

Thus, the study of the interconnections between the measure $d\tau$ and OPUC
is equivalent in a certain sense to the study of the interconnections between the measure and  {\it Verblunsky coefficients}
or Szeg\"o recurrences. 

It is important that the Toeplitz matrices (or measures) generate not only Szeg\"o recurrences but also discrete 
Dirac systems, which are dual to Szeg\"o recurrences \cite{FKRS08}. Moreover, we show that these {\it Dirac systems
are determined by Verblunsky-type coefficients}. Therefore, the study of discrete Dirac systems and corresponding
Verblunsky-type coefficients could bring  interesting results in the theory of OPUC and vice versa.

The most essential example of such interaction one finds in the continuous case, where the seminal
paper ``Continuous analogues of propositions on polynomials orthogonal on the unit circle" by M.G. Krein
\cite{Krein} led to solving of the inverse spectral problem for continuous Dirac-type systems. 
Here Krein system may be considered as an analog of the Szeg\"o recurrence and both inverse
spectral problems (for Krein system in \cite{Krein} and for Dirac-type systems in \cite{ALS02, 
SaL3}) are solved using continuous analogues of Toeplitz matrices. On the inverse spectral problems for
continuous Dirac-type systems see also the related results and references in \cite{EGRST, ALSJST, SaSaR}.

Another example of the mentioned above interaction is connected with the interesting paper \cite{DerSi}, where the nonclassical (indefinite) analog
of Szeg\"o theorem from our work \cite{ALSJFA} is derived (for the scalar case) using the theory of orthogonal polynomials.
We note that the nonclassical  analog of Szeg\"o theorem is derived in \cite{ALSJFA} for the case of block Toeplitz matrices
(and not only for the scalar case of Toeplitz matrices). It would be useful to obtain this result  in full generality 
via the theory of  orthogonal polynomials (namely, {\it matrix} orthogonal polynomials) as well.

Similar to \cite{ALSJFA}, we consider here the case of block Toeplitz matrices $S(n)$
 (and the Hankel matrices $H(n)$ are  also block Hankel matrices). Correspondingly, the functions $\tau(t)$ are matrix functions. 

Since the analogies between Szeg\"o recurrences and discrete Dirac  systems are of interest, in Section \ref{Verb1} we present
 Verblunsky-type (instead of Verblunsky) matrix coefficients (see Definition \ref{DnVerb1})
and  prove {\it Verblunsky-type
results} for Dirac systems and Toeplitz matrices (see Theorems \ref{Tm1}, \ref{Tm2} and \ref{TmDV}).
The related procedure to recover Dirac system from the Weyl function is formulated in Theorem \ref{TmWD}.

In the case of the orthogonal polynomials on the real axis (OPRL), Jacobi matrices appear instead of the
Szeg\"o recurrences and certain canonical systems appear instead of the Dirac systems.
In Section 3, we introduce the corresponding Verblunsky-type coefficients (see Definition \ref{DnVerb2})
and prove Verblunsky-type results for canonical systems and Hankel matrices (see Theorem \ref{TmM2},
Corollary \ref{CyL} and Remark \ref{RkL}). The related results on spectral functions of the canonical systems
are formulated in Theorem \ref{TmF}. 

The appendix contains a proof of the important for the theory
interpolation Theorem \ref{Hamburg}.

We note that the one to one mappings between Dirac systems and Szeg\"o recurrences are given in
\cite[Theorem 2.5]{FKRS14}. The one to one correspondence between canonical systems \eqref{Ve36}, \eqref{Ve37} 
and block Jacobi matrices is discussed in \cite[Section 8.2]{SaL2}. However, the correspondence between
Verblunsky and Verblunsky-type coefficients is rather complicated. It is better to choose which
coefficients to use depending on the problem.

As usual $\BC$ stands for the complex plane, $\BC_+$ ($\BC_-$) stands for the upper (lower)
open complex half-plane, and $\BC^{p\times r}$ stands for the set of $p \times r$ matrices with
complex-valued entries. We  often use the same notations in Sections \ref{Verb1} and \ref{Verb2} 
for the notions which play similar roles
in the schemes  for  Toeplitz and Hankel matrices, respectively.

\section{Verblunsky-type
theorems \\ for Dirac systems}\label{Verb1}
 \setcounter{equation}{0}
 \subsection{Toeplitz matrices and Dirac systems}\label{subD}
 Let $\tau(t)$ be an nondecreasing $p\times p$ matrix function on $[-\pi, \, \pi]$ such that the block Toeplitz matrix
\begin{align}\label{c0}
S(n)=\{s_{j-i}\}_{i,j=1}^n, \quad s_k:=\frac{1}{2 \pi}\int_{-\pi}^{\pi}\E^{\I k t}d\tau(t)
\end{align}
is positive-definite (i.e. $S(n)>0$). For any $S(n)=S(n)^*$ the following matrix identity
is valid (see \cite{ALSJFA} and references therein):
\begin{align}& \label{c1}
AS(n)-S(n)A^*= \I \Pi J \Pi^*;  \quad \Pi=\begin{bmatrix}\Phi_1 & \Phi_2\end{bmatrix},  
\\  & \label{c2}
A= \left\{ a_{j-i}^{\,} \right\}_{i,j=0}^n,
          \quad  a_k  =  \left\{ \begin{array}{lll}
                                  0 \, & \mbox{ for }&
k> 0  \\
\displaystyle{\frac{\I}{ {\, 2 \,}}} \,
                                I_p
                                  & \mbox{ for }& k =
0  \\
                                \, \I \, I_p
                                  & \mbox{ for }& k <
0
                          \end{array} \right.
, \qquad J=\begin{bmatrix} 0 & I_p \\ I_p & 0 \end{bmatrix};
\\ & \label{c3}
\Phi_1 = \left[
\begin{array}{c}
I_{p}  \\ I_p \\ \cdots \\ I_{p}
\end{array}
\right], \quad \Phi_2 =  \left[
\begin{array}{l}
s_0/2  \\ s_0/2 + s_{-1} \\ \cdots \\ s_0/2+ s_{-1} + \ldots +
s_{1-n}
\end{array}
\right] +\I \Phi_1 \nu, \quad \nu = \nu^*;
\\ & \label{c3'}
A=A(n), \quad \Pi=\Pi(n), \quad \Phi_1=\Phi_1(n), \quad \Phi_2=\Phi_2(n),
\end{align}
where $\I$ is  the imaginary unit and $I_p$ is the $p\times p$ identity matrix. 
We omit the variable $n$ in $S(n), \, A(n)\in \BC^{pn \times pn}$ and in some other notations (and write simply $S$ or $A$) when the order of the matrix is clear
from the context.

\begin{Rk} One can see from \eqref{c3} that the $np \times p$ matrix $\Phi_2(n)$ is not determined by $S(n)$ uniquely
but up to a self-adjoint $p \times p$ matrix $\nu$. This $\nu$ is essential and further $($in Theorem \ref{Tm1}$)$ we establish
one to one correspondence between the pairs $\{S(n), \, \nu\}$ and Dirac systems \eqref{c10}, \eqref{c11}. If we fix
$\nu$ $($e.g., set $\nu=0)$, we should speak in Theorem \ref{Tm1} about Dirac systems with the {\it potentials} $\{C_k\}$, 
where $C_0$ has 
a special form.  
\end{Rk}

The transfer matrix function $w_A$ in Lev Sakhnovich form \cite{SaL1}  
is given, for the case of the Toeplitz matrix $S(n)$ and the identity \eqref{c1}, by the formula:
\begin{align}& \label{c4}
w_A(n, \la)= I_{2p}-\I J\Pi(n)^*S(n)^{-1}\big(A(n)-\la I_{np}\big)^{-1}\Pi(n).
\end{align} 
Since $S(n)>0$, we have $S(k)>0$ $(1 \leq k \leq n)$, and so all the matrices $S(k)$
are invertible and 
\begin{align}& \label{c8}
t_k:=\begin{bmatrix} 0  & \ldots & 0 & I_p \end{bmatrix}S(k)^{-1}\begin{bmatrix} 0  & \ldots & 0 & I_p \end{bmatrix}^*>0.
\end{align}
Therefore, according to the general factorization theorem 
for transfer matrix functions $w_A$ \cite{SaL1} (see also \cite[Theorem 1.16]{SaSaR} and further references therein),
we have factorization
\begin{align}& \label{c5}
w_A(n,\la)=w_n(\la)w_{n-1}(\la) \ldots w_1(\la) \quad (n \in \BN), 
\end{align}
where
\begin{align}& \label{c6}
w_k(\la):=I_{2p}-\I \left(\frac{\I}{2}-\la\right)^{-1}J\begin{bmatrix} X_k^* \\ Y_k^* \end{bmatrix} t_k^{-1}\begin{bmatrix} X_k & Y_k  \end{bmatrix} ,
\\ & \label{c7}
\begin{bmatrix} X_k & Y_k  \end{bmatrix}=\begin{bmatrix} 0  & \ldots & 0 & I_p \end{bmatrix}S(k)^{-1}\begin{bmatrix} \Phi_1(k) & \Phi_2(k)  \end{bmatrix}.
\end{align}
The factorization above (for the case of Toeplitz matrices) was considered and applied in \cite{FKRS08, ALSJFA} (see \cite[p. 468]{ALSJFA} and \cite[p. 219]{FKRS08}).

Next, introduce $W_{k}$ by the equality (compare with \cite[(5.42)]{FKRS08}):
\begin{align}& \label{c9}
W_{k}(\la)=\la^{-k}(\la +\I)^{k}K^*w_A(k,-\la/2)K \quad (0<k \leq n), 
\\ & \label{c9'} W_0(\la):=I_{2p},
\quad K: \frac{1}{\sqrt{2}}\begin{bmatrix}I_p & -I_p \\ I_p & I_p \end{bmatrix}.
\end{align}
Using \eqref{c5}-\eqref{c9'} we will prove the following proposition.
\begin{Pn}\label{PnDir} The matrix function $W_{k}(\la)$ given by \eqref{c9} $($and generated by the Toeplitz matrix $S(n)>0)$
is a fundamental solution of the discrete Dirac system
\begin{align}& \label{c10}
W_{k+1}(\la)-W_k(\la)=-\frac{\I}{\la}j C_kW_k(\la), \quad j:=\begin{bmatrix}I_p & 0 \\ 0 & -I_p \end{bmatrix};
\\ & \label{c11}
C_k>0, \quad C_k j C_k=j \quad (0 \leq k<n).
\end{align}
Moreover, we have
\begin{align}& \label{c12}
C_k=2K^*\b(k)^*\b(k)K-j, \quad \b(k):=t_{k+1}^{-1/2}\begin{bmatrix} X_{k+1} & Y_{k+1}  \end{bmatrix}.
\end{align}
\end{Pn}
\begin{proof}. It easy to see that
\begin{align}& \label{c13}
K^*=K^{-1}, \quad K^*J K=j.
\end{align}
Hence, in view of \eqref{c5}, \eqref{c9} and \eqref{c9'}, we have
\begin{align}& \label{c14}
W_{k+1}(\la)=\frac{\la +\I}{\la}K^*w_{k+1}(-\la/2) K W_k(\la) \quad (k \geq 0).
\end{align}
Taking into account \eqref{c6}, \eqref{c12}  and \eqref{c13}, we derive
\begin{align} \nn
\frac{\la +\I}{\la}K^*w_{k+1}(-\la/2)K &=\frac{\la +\I}{\la} I_{2p}-\frac{2\I}{\la}K^*JK K^*\begin{bmatrix} X_{k+1}^* \\ Y_{k+1}^* \end{bmatrix} t_{k+1}^{-1}
\begin{bmatrix} X_{k+1} & Y_{k+1}  \end{bmatrix}K
\\ & \label{c15}
=I_{2p}-\frac{\I}{\la}j C_k.
\end{align}
Thus, according to \eqref{c14} and \eqref{c15}, the relations \eqref{c10} and \eqref{c12} are valid, and it remains to prove
\eqref{c11} for $C_k$ given by \eqref{c12}.

In view of the particular case (i.e., the case $S(k+1)=S(k+1)^*$) of \cite[Proposition 2.1]{ALSJFA}, we have
$t_{k+1}=X_{k+1}Y_{k+1}^*+Y_{k+1}X_{k+1}^*$, that is,
\begin{align}& \label{c16}
\b(k)J\b(k)^*=I_p
\end{align}
for $\b(k)$ introduced in \eqref{c12}. (The fact follows after some transformations from the identity \eqref{c1}.)
Using \eqref{c13}, we rewrite \eqref{c12} and \eqref{c16} in the form
\begin{align}& \label{c17}
C_k=2\wh \b(k)^* \wh \b(k)-j, \quad \wh \b(k) j \wh\b(k)^*=I_p, \quad \wh \b(k):= \b(k) K.
\end{align}
The second equality in \eqref{c17} yields
\begin{align}& \label{c18}
\wh \b(k)^* \wh \b(k)j \wh \b(k)^* \wh \b(k)=\wh \b(k)^* \wh \b(k).
\end{align}
According to the first equality in \eqref{c17} and to \eqref{c18}, the relation $C_k j C_k=j$ (i.e., the second relation in \eqref{c11})
holds.

Finally, in order to prove $C_k>0$ (i.e.,  the first relation in \eqref{c11}) we use some $p \times 2p$ matrix $\breve \b(k)$
such that
\begin{align}& \label{c19}
\wh \b(k) j \breve \b(k)^*=0, \quad \breve \b(k) j \breve \b(k)^*=- I_p.
\end{align}
The existence of $\breve \b(k)$ easily follows from \cite[Proposition 1.44]{SaSaR}. Taking into account \eqref{c19} and
$\wh \b(k) j \wh\b(k)^*=I_p$, we derive
\begin{align}& \label{c20}
\begin{bmatrix} \wh \b(k) \\ \breve \b(k) \end{bmatrix} j \begin{bmatrix} \wh \b(k) \\ \breve \b(k) \end{bmatrix}^*=
j=\begin{bmatrix} \wh \b(k) \\ \breve \b(k) \end{bmatrix}^* j \begin{bmatrix} \wh \b(k) \\ \breve \b(k) \end{bmatrix}.
\end{align}
Let us substitute the right-hand side of \eqref{c20} for $j$ in the first equality in \eqref{c17}. We obtain
\begin{align}\nn
C_k&=2\wh \b(k)^* \wh \b(k)-\begin{bmatrix} \wh \b(k) \\ \breve \b(k) \end{bmatrix}^* j \begin{bmatrix} \wh \b(k) \\ \breve \b(k) \end{bmatrix}
\\ & \label{c21}
=\wh \b(k)^* \wh \b(k) +\breve \b(k)^* \breve \b(k)=\begin{bmatrix} \wh \b(k) \\ \breve \b(k) \end{bmatrix}^* \begin{bmatrix} \wh \b(k) \\ \breve \b(k) \end{bmatrix}>0.
\end{align}
\end{proof}
The inverse to Proposition \ref{PnDir} statement, namely, the statement that each discrete Dirac system \eqref{c10}, \eqref{c11} is generated by some
block Toeplitz matrix $S(n)>0$ (and the procedure to recover such $S(n)$) easily follows from \cite[Theorem 5.2]{FKRS08}.
\begin{Pn}\label{PnRT} For each Dirac system \eqref{c10}, \eqref{c11} there exist a block Toeplitz matrix $S(n)>0$ and a $p \times p$ matrix $\nu=\nu^*$
$($which is used in the definition \eqref{c4} of $\Phi_2(n))$
such that relations \eqref{c12} hold.
\end{Pn}
\begin{proof}. According to \cite[Theorem 5.2]{FKRS08}, there are $S(n)>0$ and $\nu = \nu^*$ such that the corresponding $C_0$ is given, indeed,
by the formula \eqref{c12}, and the matrices $C_k$ $(0 <k<n)$ are given by the formulas
\begin{align}\label{c22}&
C_k=2K^*\vt(k) K-j  \quad (k>0), \\
\label{c23}&
 \vt(k):=\Pi(k+1)^*S(k+1)^{-1}\Pi(k+1)-\Pi(k)^*S(k)^{-1}\Pi(k).
\end{align}
We note that the notations in \cite{FKRS08} slightly differ from the notations in this paper
and, in particular, $P_kS P_k^*$ in \cite{FKRS08} is denoted by $S(k+1)$ here.
Moreover, for $\a_0$ from \cite[Theorem 5.2]{FKRS08} we have $\a_0=\frac{s_0}{2}+\I \nu$.

Writing $S(k+1)$ in the block form
\begin{align}& \label{c24}
S(k+1)=\begin{bmatrix} S(k) & S_{12} \\ S_{21} & s_0 \end{bmatrix}
\end{align}
and using formula \cite[(2.7)]{ALSJFA} for $T(k+1)=S(k+1)^{-1}$, we obtain
\begin{align}& \label{c25}
\begin{bmatrix} 0  & \ldots & 0 & I_p \end{bmatrix}S(k+1)^{-1}\Pi(k+1)=t_{k+1}\begin{bmatrix} S_{21}S(k)^{-1} & -I_p \end{bmatrix}\Pi(k+1),
\end{align}
and we rewrite also \eqref{c23} in the form
\begin{align}& \label{c26}
 \vt(k)=\Pi(k+1)^* \begin{bmatrix} S(k)^{-1}S_{12} \\ -I_p \end{bmatrix}t_{k+1}\begin{bmatrix} S_{21}S(k)^{-1} & -I_p \end{bmatrix}    \Pi(k+1).
\end{align}
Relations \eqref{c25} and \eqref{c26} yield
\begin{align}& \label{c27}
\vt(k)=\Pi(k+1)^*S(k+1)^{-1} \begin{bmatrix} 0  \\ \ldots \\ 0 \\ I_p \end{bmatrix}t_{k+1}^{-1} \begin{bmatrix} 0  & \ldots & 0 & I_p \end{bmatrix}S(k+1)^{-1}\Pi(k+1).
\end{align}
On the other hand, formula \eqref{c7} and the definition of $\b(k)$ in \eqref{c12} imply that
\begin{align}& \label{c28}
\b(k)=t_{k+1}^{-1/2}\begin{bmatrix} 0  & \ldots & 0 & I_p \end{bmatrix}S(k+1)^{-1}\Pi(k+1)
\end{align}
From \eqref{c27} and \eqref{c28}, it is immediate that $\vt(k)=\b(k)^*\b(k)$, and so formula
\eqref{c22} coincides with \eqref{c12} (for $k>0$). It means that relations \eqref{c12} are valid, and $S(n)$ constructed in \cite[Theorem 5.2]{FKRS08}
generates (in the sense of Proposition \ref{PnDir}) the initial Dirac system.
\end{proof}
\begin{Rk}\label{Uniq} It is essential that only one pair $\{S(n)>0, \,\, \nu=\nu^*\}$ generates $($in the sense of Proposition \ref{PnDir}$)$ each Dirac system
\eqref{c10}, \eqref{c11}. Indeed, according to \cite[Theorem 4.1]{ALSJFA} or \cite[Theorem 6.1]{FKRS08} each $w_A(n,\la)$ uniquely determine
$S(n)$ $($and $\a_0=\frac{s_0}{2} +\I \nu)$, which generate it, and $w_A(n,\la)$ is in turn uniquely recovered from the Dirac system
\eqref{c10}, \eqref{c11} using \eqref{c9}.
\end{Rk}

The following theorem is immediate from Propositions \ref{PnDir}, \ref{PnRT} and Remark \ref{Uniq}.
\begin{Tm} \label{Tm1} For each $n \in \BN$, there is a one to one correspondence between the pairs $\{S(n), \, \nu\}$
$($of  the block Toeplitz matrices $S(n)>0$ and the arbitrary $p \times p$ matrices $\nu=\nu^*$
which are used in the definition \eqref{c4} of $\Phi_2(n))$
and Dirac systems \eqref{c10}, \eqref{c11}. This correspondence is given by the formula \eqref{c12}.
\end{Tm}
\subsection{Verblunsky-type coefficients and Dirac systems}\label{subV}
\paragraph{1.} Further in the text, we introduce Verblunsky-type coefficients $\rho_k$ for Dirac systems via matrices $C_k$. The connection between
$\rho_k$ and $C_k$ is similar to the connection between Verblunsky coefficients $a_k$ and Szeg\"o recurrences \eqref{03}.
The next theorem is a special case of  \cite[Proposition 2.4]{FKRS14}. Note that the paper \cite{FKRS14} is dedicated
to Dirac systems with $j=\begin{bmatrix} I_{p_1} & 0 \\ 0 &  -I_{p_2} \end{bmatrix}$, and Dirac systems considered here
appear in \cite{FKRS14} when one puts $p_1=p_2=p$.
\begin{Pn}\label{novHE} Let a $2p \times 2p$ matrix $C$ be $j$-unitary and positive. Then, it admits a representation
\begin{align}\label{nov3}&
C= { \cal D}  F,
\end{align}
where   $F$ and $  {\cal D}$ are of the form 
\begin{align}\label{nov1}&
 F= \left[
\begin{array}{cc}
I_{p} &  \rho \\   \rho^* & I_{p}
\end{array}
\right],
\quad
   {\mathcal D}= {\mathrm{diag}}\Big\{
\big(
I_{p}-  \rho  \rho^* \big)^{-\frac{1}{2}}, \, \,
\big(I_{p}- \rho^* \rho\big)^{-\frac{1}{2}}\Big\}, 
\\ \label{nov2}&
 \rho^*\rho <I_{p}.
\end{align}
\end{Pn}

\begin{Rk} \label{RkHalm}
Note that $F\cld$ is a well-known expression $($such that $F\cld=\cld F)$ which is sometimes called the Halmos extension of $\rho$
$($see, e.g., \cite[p. 167]{DFK}$)$.
\end{Rk}

\begin{Rk} \label{RkUnrho}
Clearly, there is a unique $\rho$ satisfying \eqref{nov3}. Namely, partitioning $C$ into $p \times p$ blocks $C=\{c_{ij}\}_{i,j=1}^2$, we express  $\rho$ by the formula$:$
\begin{align}\label{nov4}&
\rho=c_{11}^{-1}c_{12}.
\end{align}
\end{Rk}

Using Proposition \ref{novHE} and Remark \ref{RkUnrho}, we introduce the notion of Verblunsky-type coefficients.
\begin{Dn} \label{DnVerb1}
Verblunsky-type coefficients $\rho_k$ of a block Toeplitz matrix $S(n)>0$  $($more precisely, of the pair $\{S(n), \, \nu\})$
or, equivalently,  of the corresponding Dirac system \eqref{c10}, \eqref{c11} are given, using representations \eqref{nov3} of the matrices $C_k$,
by  the formulas
\begin{align}\label{nov5}&
\rho_k=\left(\begin{bmatrix} I_p & 0 \end{bmatrix}C_k\begin{bmatrix} I_p \\ 0 \end{bmatrix}\right)^{-1}\begin{bmatrix} I_p & 0 \end{bmatrix}C_k
\begin{bmatrix} 0 \\ I_p \end{bmatrix} \qquad (0\leq k<n),
\end{align}
where the matrices $C_k$ are expressed via $S(n)$ in \eqref{c12} $($using formulas \eqref{c8} and \eqref{c7}$)$.
\end{Dn}
In view of Theorem \ref{Tm1}, Proposition \ref{novHE} and Remark \ref{RkUnrho},   we can formulate
the theorem which confirms the correctness of the Definition \ref{DnVerb1}.
\begin{Tm} \label{Tm2} The Verblunsky-type coefficients satisfy the inequalities 
$$\|\rho_k\|~<~1.$$ 
For each $n \in \BN$, there is a one to one correspondence between 
the pairs $\{S(n), \, \nu\}$
$($of  the block Toeplitz matrices $S(n)>0$ and  of the arbitrary $p \times p$ matrices $\nu=\nu^*$
which are used in the definition \eqref{c4} of $\Phi_2(n))$
 and  the sequences of the coefficients $\rho_k$ such that $\|\rho_k\|<1$ $(0\leq k<n)$. This correspondence is given by the formula \eqref{nov5}.
\end{Tm}
\paragraph{2.} Using considerations from Proposition \ref{PnDir} and Remark \ref{Uniq} we see that if the sequence of Verblunsky-type coefficients $\{\rho_k\}$ $(0 \leq k<N)$ corresponds to the pair $\{S(N)=\{s_{j-i}\}_{i,j=1}^N, \, \nu\}$ (and vice versa), then the sequence $\{\rho_k\}$ $(0 \leq k<n<N)$ corresponds to the pair $\{S(n), \, \nu\}$ with the reduced
matrix $S(n)=\{s_{j-i}\}_{i,j=1}^n$ (and the pair $\{S(n), \, \nu\}$  corresponds to this sequence). Together with Theorem \ref{Tm2} it means
that there is a one to one correspondence between the infinite sequences of Verblunsky-type coefficients $\{\rho_k\}$ $(0 \leq k<\infty)$ and the pairs consisting of $\nu=\nu^*$ and the infinite sequences
of the blocks $\{s_k\}$ $(0\leq k<\infty)$ such that $S(n)=\{s_{j-i}\}_{i,j=1}^n>0$ for all $n \geq 1$.
Moreover, in the case $S(n)\geq 0$ for all $n \geq 1$,  there is a unique measure $d \tau$ such that \eqref{c0} holds (see, e.g., \cite[Theorem 1]{DGK}).
Now, Verblunsky-type theorem below follows from Theorem \ref{Tm2}.
\begin{Tm}\label{TmDV}
There is a one to one correspondence between the pairs consisting of the measures $d\tau$, such that $S(n)> 0$ for all $S(n)$ generated by \eqref{c0} $(n \geq 1)$,
and arbitrary $p\times p$ matrices $\nu = \nu^*$ on one side, and 
the sequences $\{\rho_k\}$ $(0 \leq k<\infty)$, where $\| \rho_k \|<1$ $($Verblunsky-type coefficients of the measure$)$, on the other side.
\end{Tm}
\paragraph{3.}  The discussed above correspondence between Verblunsky-type coefficients and Dirac systems on one side and block Toeplitz matrices on the
other side is essential for solving inverse spectral problems. We demonstrate this by the result, which is immediate from
\cite[Theorem 5.2 and Section 6]{FKRS08} and from the proof of our Proposition \ref{PnRT}. We consider in this paragraph Dirac systems 
on the semiaxis $k \in \left\lbrace 0,1,2,
\ldots\right\rbrace $ $($i.e., on the semiaxis $k\geq 0)$ and substitute conditions \eqref{c11} on the finite interval by the similar
conditions on the semiaxis:
\begin{align}& \label{c11'}
C_k>0, \quad C_k j C_k=j \quad (0 \leq k<\infty).
\end{align}
\begin{Dn}\label{DnWD}
 A $p \times p$ matrix function
$\vp(\la)$
holomorphic in the lower complex half-plane $\BC_-$ is called a Weyl function
for Dirac system \eqref{c10}, \eqref{c11'}
  if the
inequality
\begin{equation} \label{c30}
\sum_{k=0}^\infty [\I \vp(\lambda)^* \quad I_p
]q(\lambda)^k K
W_k(\lambda)^*C_kW_k(\lambda)K^*\left[\begin{array}{c}
-\I \vp(\lambda) \\ I_p
\end{array}
\right]<\infty ,
\end{equation}
holds for
$q(\lambda):=|\lambda^2|(|\lambda^2|+1)^{-1}$.
\end{Dn}
\begin{Tm}\label{TmWD} There is a unique Weyl function $\vp(\la)$ of Dirac system \eqref{c10}, \eqref{c11'}.
This Weyl function admits representation
\begin{equation}\label{c31}
\I\vp\Big(\I\frac{(z+1)}{(z-1)}\Big)=\a_0+\sum_{k=1}^{\infty}s_{-k}z^k,
\end{equation}
where $\{s_{j-i}\}_{i,j=1}^n >0$  for $s_0=\a_0+\a_0^*$ and for all $0 < n<\infty$.
Each function $\vp(\la)$, such that \eqref{c31} holds and $\{s_{j-i}\}_{i,j=1}^n >0$  for $s_0=\a_0+\a_0^*$ and all $0 < n<\infty$,
is a Weyl function of some Dirac system \eqref{c10}, \eqref{c11'}.

The Dirac system \eqref{c10}, \eqref{c11'} is
uniquely recovered from
the Weyl
function via the formula \eqref{c12} using definitions \eqref{c3}, \eqref{c8}, \eqref{c7} and the equality $\nu=(\a_0-\a_0^*)/(2\I)$.

\end{Tm}
%%%%%%%%%%%%%%%%%%%%%%%%%%%%%%%%%%%%%%%%%%%%%%%%%%%%%%%%%%%%%%%%%%%
\section{Verblunsky-type theorems \\ for canonical systems}\label{Verb2} 
 \setcounter{equation}{0}
 \subsection{Hankel matrices and canonical systems}\label{subH} 
 Consider a self-adjoint $n \times n$ block Hankel matrix $H$ with the $p \times p$ 
 blocks $H_k$ $(0 \leq k \leq 2n-2)$:
\begin{align}& \label{H1}
H=H(n)=\{H_{i+j-2}\}_{i,j=1}^n \quad (n \geq 1).
\end{align}
The self-adjoint Hankel matrix (equivalently, the Hankel matrix where $H_k=H_k^*$) satisfies the matrix identity
 \begin{align}& \label{H2}
AH-HA^*=\I \Pi J \Pi^*, \quad \Pi=\Pi(n)=\begin{bmatrix}\Phi_1(n) & \Phi_2(n)\end{bmatrix},
\\ & \label{H3}
A=A(n)=\begin{bmatrix} 0 & 0 & \ldots & 0 & 0 \\ I_p & 0 & \ldots & 0 & 0 \\  \cdots & \cdots & \cdots & \cdots &\cdots 
\\ 0 & 0 &\ldots & I_p & 0
\end{bmatrix}, \quad J=\begin{bmatrix} 0 & I_p \\ I_p & 0 \end{bmatrix},
\\ & \label{H4}
\Phi_1=\Phi_1(n)=-\I \begin{bmatrix} 0 \\ H_0 \\ H_1 \\ \cdots \\ H_{n-1}\end{bmatrix}, \quad 
\Phi_2=\Phi_2(n)=\begin{bmatrix} I_p \\  0\\  0\\ \cdots \\ 0\end{bmatrix}.
\end{align}
Clearly $A$ and $H$ are $n p \times np$ matrices and $\Phi_1$ and $\Phi_2$
are $ np \times p$ matrices. When we deal with the matrices $A$, $H$ or $\Pi$ with various numbers of block rows,
we write $H(n)$, $A(n)$  and $\Pi(n) $ instead of  $H$, $A$  and $\Pi $, respectively. 

We note that we often use the same notations in Sections \ref{Verb1} and \ref{Verb2} for the notions which play similar roles
in the schemes  for  Toeplitz and Hankel matrices, respectively.
One can see that in the case of
the Toeplitz matrices $S$ we had $\Phi_1$ independent of the choice of $S$ and $\Phi_2$ depending
on this choice, but in the case of
the Hankel matrices, $\Phi_2$ does not depend on $H$ and $\Phi_1$ depends (see \eqref{H4}). This difference
(which is caused by the literature used in Sections \ref{Verb1} and \ref{Verb2}) is mathematically insignificant.

In view of \eqref{H2}, an explicit solution of the interpolation problem  (i.e., of the truncated
Hamburger moment problem) for $H>0$ is easily obtained via the method of operator identities \cite{SaL0, SaL1, SaL2, SaL3}
using V.P. Potapov's fundamental matrix inequalities
\cite{Am, Pot1}. 

First, introduce {\it several important notations}. Recall that the transfer matrix function from \cite{SaL1}  
is given (for the case of the identity \eqref{H2}) by the formula:
\begin{align}& \label{H5}
w_A(n, \la)= I_{2p}-\I J\Pi(n)^*H(n)^{-1}\big(A(n)-\la I_{np}\big)^{-1}\Pi(n).
\end{align}

A $2p \times p$ matrix-valued function
(matrix function)
$\clq(z)$, which is holomorphic in the open upper half-plane $\BC_+$, is called {\it nonsingular with property-$J$}
if
\begin{align}& \label{H6}
\clq(z)^*\clq(z)>0, \quad \clq(z)^* J \clq(z)\geq 0. 
\end{align}
The matrix function $\clq$ is often partitioned into two $p \times p$ matrix functions,
in which case one speaks about a nonsingular pair with property-$J$.

Finally, let us {\it denote by} $\cle(n)$ the class of nondecreasing $p \times p$ matrix functions
$\tau(t)$ ($t \in \BR$) such that
\begin{align}& \label{H7}
\int_{-\infty}^{\infty}t^{2n-2}d\tau(t) <\infty.
\end{align}

Now, we set
\begin{align}& \label{H7+}
\mfa(n, z)=w_A(n,1/\ov{z})^*=w_A(1/\ov{z})^*,
\end{align}
and collect the results from \cite[pp. 19, 20]{SaL3} into the following theorem.
\begin{Tm}\label{Hamburg} Assume that the block Hankel matrix $H=H(n)$ is positive-definite
$($i.e., $H>0)$. 

Then, the matrix functions $\vp(z)=\vp_{\clq}(z)$ given by the linear fractional
transformations
\begin{align}& \label{H8}
\vp(z)=\I \mfa_1(z)\clq(z)\big(\mfa_2(z)\clq(z)\big)^{-1}, \\
& \label{H9}
 \mfa_1(z):=\begin{bmatrix}  I_p & 0 \end{bmatrix}w_A(1/\ov{z})^*,
\quad \mfa_2(z):=\begin{bmatrix}  0 & I_p \end{bmatrix}w_A(1/\ov{z})^*,
\end{align}
where the matrix functions $\clq$ are nonsingular with property-$J$ and the transfer matrix function $w_A$ is 
introduced in \eqref{H5}, belong to the Herglotz class, that is, $\Im \big(\vp(z)\big) \geq 0$ for $z\in \BC_+$.
Moreover, each matrix function $\vp(z)$ admits a unique Herglotz representation of the form
\begin{align}& \label{H10}
\vp(z)=\int_{-\infty}^{\infty}(t-z)^{-1}d\tau(t) <\infty,
\end{align}
and $\tau$ in this representation belongs $\cle(n)$.

These and only these  matrix functions $\tau \in \cle(n) $ $($i.e., $\tau$ given by \eqref{H8}--\eqref{H10}$)$ satisfy the equalities
\begin{align}& \label{H11}
H_k=\int_{-\infty}^{\infty}t^{k}d\tau(t) \quad (0 \leq k <2n-2),
\end{align}
and the inequality 
\begin{align}& \label{H12}
H_{2n-2} \geq \int_{-\infty}^{\infty}t^{2n-2}d\tau(t). \end{align}
\end{Tm}
The proof of the theorem above is based on various statements from several sections in \cite[Ch. 1]{SaL2},
and we explain in the Appendix which statements from \cite[Ch. 1]{SaL2} are used (without
formulating some of them directly).

We partition  self-adjoint  matrices $H(n)$ $(n>1)$ of the form \eqref{H1} into the following blocks
\begin{align}& \label{Hf0-}
H(n)=\begin{bmatrix} H(n-1) & \clh_{12}(n) \\
\clh_{21}(n) &  H_{2n-2} \end{bmatrix}, \quad \clh_{12}(n):=\begin{bmatrix} H_{n-1} \\ \cdots \\ 
  H_{2n-3} \end{bmatrix}, \quad \clh_{21}(n) =\clh_{12}(n) ^*.
\end{align}
Assuming
\begin{align}& \label{Hf0}
\det H(n)\not=0, \quad \det H(n-1)\not=0,
\end{align}
we put
\begin{align}& \label{Hf0+}
T(m)= H(m)^{-1},  \quad \clh_{12}=\clh_{12}(n), \quad \clh_{21}=\clh_{21}(n), \\
& \label{Hf0++}
  t_n:=\big(H_{2n-2}-\clh_{21}T(n-1)\clh_{12}\big)^{-1}.
\end{align}
Then, the inversion formula below easily follows from the representation \eqref{Hf0-} (and is given in \cite[(2.7)]{ALSJFA}):
\begin{align}& \label{Hf1}
T(n)=\begin{bmatrix} T(n-1)+T(n-1)\clh_{12}t_n\clh_{21}T(n-1) & -T(n-1)\clh_{12}t_n \\
-t_n\clh_{21}T(n-1) & t_n \end{bmatrix}.
\end{align}

Applying general factorization theorem for transfer matrix functions $w_A$ \cite{SaL1} (see also \cite[Theorem 1.16]{SaSaR}) 
to the case of the Hankel matrix $H(n)$ and using \eqref{H3}, \eqref{Hf0+} and \eqref{Hf1},
we have the factorization 
\begin{align}& \label{Hf2}
w_A(n,\la)=\clw(n,\la)w_A(n-1, \la), 
\\ & \label{Hf3}  \clw(n,\la)=I_{2p}+\frac{\I}{\la}J \Pi(n)^*T(n)P_2^*t_n^{-1}P_2T(n)\Pi(n), 
\\ & \label{Hf3+}
P_2=P_2(n):=
\begin{bmatrix} 0 & \ldots & 0 & I_p
 \end{bmatrix} \in \BC^{p\times pn}.
\end{align}
We note that  we apply above   \cite[Theorem 1.16]{SaSaR} to the so called symmetric $S$-node (see \cite[Definition 1.12]{SaSaR}).
Similar to some other notations, we omit the variable $n$ in $P_2(n)$ and write simply $P_2$ when the order of $P_2$ is clear
from the context.

When $H(n)>0$,  the equalities $\det H(k)\not=0$ hold for $1 \leq k \leq n$, and so we may factorize successively
$w_A(n-1)$, $w_A(n-2)$, $\ldots, \, w_A(2)$.  Thus, we can rewrite \eqref{Hf2} in the form
\begin{align}& \label{Hf4}
w_A(n,\la)=w_n(\la)w_{n-1}(\la) \ldots w_1(\la), \\
& \label{Hf4+}
 w_k(\la):=\clw(k,\la) \quad (k>1), \quad w_1(\la):=w_A(1, \la).
\end{align}
It is  also easy to see that  in the case $H(n)>0$ we have
\begin{align}& \label{Hf5}
t_k>0 \quad (1 \leq k \leq n), \quad w_A(1, \la)=I_{2p}+\frac{\I}{\la}J\begin{bmatrix} 0 \\ I_p
 \end{bmatrix} H_0^{-1}\begin{bmatrix} 0 & I_p
 \end{bmatrix} .
\end{align}
It follows from \eqref{Hf3}--\eqref{Hf5} that the system
\begin{align}& \label{Hf6}
y_{k+1}=w_{k+1}(\la)y_{k} \qquad (0 \leq k < n)
\end{align}
is a discrete canonical system (which is generated by the Hankel matrix $H(n)$). 
More precisely, we have the following proposition.
\begin{Pn}\label{PnCan} Each Hankel matrix $H(n)>0$ determines a canonical system
\begin{align}& \label{Ve1}
y_{k+1}=w_{k+1}(\la)y_{k} \quad (0 \leq k< n),  \quad w_{k+1}( \la)=I_{2p}+\frac{\I}{\la}J Q_{k};
\\ & \label{Ve2} 
 Q_k=\Pi(k+1)^*T(k+1)P_2^*t_{k+1}^{-1}P_2T(k+1)\Pi(k+1) \quad (1 \leq k),
\\ & \label{Ve2+} 
Q_0=
\begin{bmatrix} 0 \\ I_p
 \end{bmatrix} H_0^{-1}\begin{bmatrix} 0 & I_p
 \end{bmatrix}, 
\end{align}
where $T(m)=H(m)^{-1}$. The  matrices $Q_k$  $(0\leq k< r<n)$, that is the first $r$ values of the Hamiltonian $\{Q_k\}$,
are determined by the blocks $H_k$ $(0\leq k\leq 2r-2)$ of $H(n)$. 
\end{Pn}
%%%%%%%%%%%%%%%%%%%%%%%%%%%%%%%%%%%%%%%%%%%%%%%%%%%%%%%%%%

 \subsection{Verblunsky-type coefficients and canonical systems}\label{subVH}
 
\paragraph{1.} According to \eqref{Ve2} and \eqref{Ve2+}, we have
\begin{align}& \label{Ve3}
Q_k= \om_k^* t_{k+1}^{-1}\om_k, \\
& \label{Ve3!}
  \om_k: =P_2 T(k+1)\Pi(k+1) \quad (0 \leq k<n).
\end{align}
Here we take into account that in the special case $k=0$ we have
\begin{align}& \label{Ve3+}
A(1)=0, \quad \Pi(1)=\begin{bmatrix} 0 & I_p
 \end{bmatrix}, \quad P_2(1)=I_p, \quad T(1)=t_1=H_0^{-1}.
\end{align}
In order to show that $\om_k$ may be considered as Verblunsky-type coefficients,
we need to study their properties. In particular, similar to the case of Toeplitz
matrices (see the proof of \eqref{c16}), we need an expression for $t_{k+1}$. 
\begin{Pn}\label{PnProp} The matrices $\om_k$ introduced in \eqref{Ve3!}
have the following properties$:$
\begin{align} & \label{Ve13}
\om_0\begin{bmatrix} I_p \\ 0
 \end{bmatrix} =0, \quad \om_0\begin{bmatrix} 0 \\ I_p 
 \end{bmatrix} =t_1; 
 \\ & \label{Ve13+}
 \om_kJ \om_k^*=0, \quad \I \om_{k}J \om_{k-1}^*=t_{k+1} \quad (0 < k <n).
\end{align}
\end{Pn}
\begin{proof}. The properties \eqref{Ve13} of $\om_0$ are immediate from \eqref{Ve3!} and \eqref{Ve3+}.

The matrix identity
\eqref{H2} implies that
\begin{align}\label{Ve7}&
T(k+1)A-A^*T(k+1)=\I T(k+1)\Pi(k+1)J \Pi(k+1)^* T(k+1).
\end{align}
In view of \eqref{H3} and \eqref{Ve7}, we obtain 
$$P_2T(k+1)\Pi(k+1)J \Pi(k+1)^* T(k+1)P_2^*=0.$$
Now, the first equality in \eqref{Ve13+} follows from the definition of $\om_k$ in \eqref{Ve3!}.

According to \eqref{H4} and the definition of $\om_{k-1}$, we have the formula
\begin{align}\label{Ve5+}&
\om_{k-1}^*=\Pi(k+1)^* \begin{bmatrix} T(k) \\ 0
 \end{bmatrix} P_2(k)^*.
\end{align}
Taking into account \eqref{Ve3!}, \eqref{Ve5+} and \eqref{H2}, 
we derive
\begin{align}\label{Ve9}
\om_k J \om_{k-1}^*
& = P_2 T(k+1)\Pi(k+1)J \Pi(k+1)^*\begin{bmatrix} T(k) \\ 0
 \end{bmatrix} P_2(k)^*
 \\ \nn &= -\I P_2 \big(T(k+1)A(k+1)H(k+1)-A(k+1)^*\big)\begin{bmatrix} T(k) \\ 0
 \end{bmatrix} P_2(k)^*.
\end{align}
Definitions \eqref{H3} and
\eqref{Hf3+} imply $P_2A^*=0$, and so we rewrite \eqref{Ve9} in the form
\begin{align}
 \label{Ve10}& 
\om_k J \om_{k-1}^*= -\I P_2(k+1) T(k+1)A(k+1)H(k+1)\begin{bmatrix} T(k) \\ 0
 \end{bmatrix} P_2(k)^* .
\end{align}
Using the partitioning \eqref{Hf0-} of $H(k+1)$ and recalling again the definitions \eqref{H3} and
\eqref{Hf3+},  we obtain
\begin{align} &
 \label{Ve11} 
A(k+1)H(k+1)\begin{bmatrix} T(k) \\ 0
 \end{bmatrix} P_2(k)^* =P_2(k+1)^*.
\end{align}
Substitute \eqref{Ve11} into \eqref{Ve10} and use \eqref{Hf1}  in order to derive
\begin{align} & \label{Ve12}
\om_k J \om_{k-1}^*=-\I P_2(k+1) T(k+1) P_2(k+1)^* =-\I t_{k+1}.
\end{align}
The second equality in \eqref{Ve13+} is immediate from  \eqref{Ve12}.
\end{proof}
Formulas \eqref{Ve3}, \eqref{Ve13} and \eqref{Ve13+} yield the corollary
\begin{Cy}\label{CyQ} The Hamiltonian $\{Q_k\}$ is expressed via $\{\om_k\}$ by the relation
\begin{align} &
 \label{Ve8} 
Q_0=\om_0^* \left(\om_0\begin{bmatrix} 0 \\ I_p 
 \end{bmatrix}\right)^{-1}\om_0, \quad Q_k=\om_k^*\big(\I \om_k J \om_{k-1}^*\big)^{-1} \om_k
 \quad (0<k<n).
\end{align} 
\end{Cy}
\begin{Dn}\label{DnVerb2} 
Verblunsky-type coefficients $\, \om_k \,$ of a block Hankel matrix $H(n)>0$  
or, equivalently,  of the corresponding canonical system \eqref{Ve1}--\eqref{Ve2+} are given
by  the formula \eqref{Ve3!}.
\end{Dn}
We note that Verblunsky-type coefficients introduced above may be considered as a parametrization
of Hankel matrices. Such parametrization is directly connected with the spectral theory and is essentially simpler than
the important Catalan or canonical Hankel parametrizations (see, e.g., \cite{Ai, FKM}).
\begin{Pn}\label{PnUniqC}
Each block Hankel matrix $H(n)>0$ is uniquely defined by its Verblunsky-type coefficients $\om_k$ $(0 \leq k <n)$.
\end{Pn}
\begin{proof}. Since $\{\om_k\}$ $(0 \leq k <n)$ are Verblunsky-type coefficients of $H(n)$, it is easy to see that each  subsequence 
 $\{\om_k\}$ $(0 \leq k <r<n)$ is the sequence of Verblunsky-type coefficients of $H(r)$. Now, we prove the proposition by induction. It follows from \eqref{Ve3+} and \eqref{Ve13}
that 
\begin{align} &
 \label{H0} 
H_0=\left(\om_0\begin{bmatrix} 0 \\ I_p 
 \end{bmatrix}\right)^{-1}. 
\end{align} 
Assume that $H(r)$ is uniquely determined by the subsequence 
 $\{\om_k\}$ $(0 \leq k <r<n)$ and let us show that the sequence $\{\om_k\}$ $(0 \leq k <r+1)$
 determines further blocks $H_{2r-1}$ and $H_{2r}$ of $H(r+1)$.
 
 Setting $n=r+1$  in formulas \eqref{H7+}, \eqref{Hf4}, \eqref{Ve1}  and in Theorem \ref{Hamburg}
 and taking into account \eqref{Ve8}, we see that $\mfa(r+1,z)$ is uniquely determined by the sequence $\{\om_k\}$ $(0 \leq k <r+1)$,
 and so the set of functions $\vp(z)$ of the form \eqref{H8} is uniquely determined by this sequence as well.
 Moreover, formulas \eqref{H10}, \eqref{H11} and representation 
$$\frac{1}{u-z}=-\frac{1}{z}\left((u/z)^{2r+1}\left(1-\frac{u}{z}\right)^{-1}+\sum_{k=0}^{2r}(u/z)^k\right)$$
yield the asymptotic formula
\begin{align} & \label{Ve14}
\vp(z)=-\sum_{k=0}^{2r-1}\frac{1}{z^{k+1}}H_k+\frac{1}{z^{2r+1}}\breve H_{2r} + o\left(\frac{1}{z^{2r+1}}\right)
\end{align}
for $z$ tending nontangentially to infinity in $\BC_+$. Thus, $\vp(z)$ and so $\{\om_k\}$ $(0 \leq k <r+1)$
determine uniquely the block $H_{2r-1}$.

Finally, given the blocks $H_0, \ldots, H_{2r-1}$ we recover  $H_{2r}$ using the second equalities in \eqref{Hf0++} and \eqref{Ve13+}:
\begin{align} &
 \label{Ve15} 
H_{2r}=(\I \om_{r}J \om_{r-1}^*)^{-1}+\clh_{21}(r+1)T(r)\clh_{12}(r+1).
\end{align} 
\end{proof}
\paragraph{2.} Some modification of the proof of Proposition \ref{PnUniqC} will show that
each sequence $\{\om_k\}$ with the properties discussed in Proposition \ref{PnProp}
is a sequence of Verblunsky-type coefficients of some Hankel matrix.
Namely, we will prove the following proposition.
\begin{Pn}\label{PnOtO} Each sequence $\{\om_k\}$ $(0 \leq k <n)$ of $p \times 2p$ matrices, such that
\begin{align} & \label{Ve16}
\om_0\begin{bmatrix} I_p \\ 0
 \end{bmatrix} =0, \quad \om_0\begin{bmatrix} 0 \\ I_p 
 \end{bmatrix} >0; 
 \\ & \label{Ve16+}
 \om_kJ \om_k^*=0, \quad \I \om_{k}J \om_{k-1}^*>0 \quad (0 < k <n),
\end{align}
is a sequence of Verblunsky-type coefficients of some block Hankel matrix $H(n)$.
\end{Pn}
\begin{proof}. Similar to the proof of Proposition \ref{PnProp} we prove this proposition by induction.
Setting $H_0=\left(\om_0\begin{bmatrix} 0 \\ I_p 
 \end{bmatrix} \right)^{-1}$, we see from \eqref{Ve3!} and \eqref{Ve3+} that  $\om_0$ is the Verblunsky-type coefficient
 of $H_0>0$. 
 
 Assume that the subsequence 
 $\{\om_k\}$ $(0 \leq k <r<n)$ is the sequence of Verblunsky-type coefficients of a block Hankel matrix $H(r)>0$.
 In order to show that $\{\om_k\}$ $(0 \leq k <r+1)$ is again a sequence of Verblunsky-type coefficients,
 we will determine the blocks $H_{2r-1}$ and $H_{2r}$ of some extension $H(r+1)$ of $H(r)$ and will prove further that 
 $\{\om_k\}$ $(0 \leq k <r+1)$
 are Verblunsky-type coefficients of this $H(r+1)$. For this, we need some preparations.
 
 Setting  $n=r$  in Theorem \ref{Hamburg}
 and putting
\begin{align} &
 \label{Ve17} 
\clq(z)\equiv \om_r^*,
\end{align}  
 we see that $\clq$ satisfies \eqref{H6} and the conditions of Theorem \ref{Hamburg} are fulfilled. Moreover, 
 $\vp(z)$ given by \eqref{H8}, where $\clq \equiv const$, is a rational matrix function and  may be considered as a rational meromorphic
 function on $\BC$.  In order to show that
\begin{align} &
 \label{Ve18} 
\vp(z)=\vp(\ov{z})^* \quad {\mathrm{for}} \quad z\in \BC, \quad \clq\equiv \om_r^*,
\end{align}  
we write a useful more general formula for $\vp$ (i.e., for $\vp_{\clq}$) and another function $\wh \vp=\vp_{\wh \clq}$  given by \eqref{H8} after substitution there of some $\wh Q$ instead of $Q$
(and of $r$ instead of $n$).
Namely, we note that in view of \cite[(1.84)]{SaSaR} we have $w_A(z)Jw_A(\ov{z})^*=J$, and so \eqref{H7+} yields
$\mfa(\ov{z})^*J\mfa(z)=J$. Hence, \eqref{H8} yields the relation
\begin{align} 
 \nn
\vp(z)-\wh \vp(\ov{z})^*&=\I\big((\mfa_2(\ov{z})\wh \clq)^{-1}\big)^*\wh \clq^*\big(\mfa_2(\ov{z})^*\mfa_1(z)+   \mfa_1(\ov{z})^*\mfa_2(z) \big)\clq  (\mfa_2(z)\clq)^{-1}
\\ \nn &=\I\big((\mfa_2(\ov{z})\wh \clq)^{-1}\big)^*\wh \clq^*\mfa(\ov{z})^*J\mfa(z)\clq  (\mfa_2(z)\clq)^{-1}
\\ \label{Ve19}  &=
\I \big((\mfa_2(\ov{z})\wh \clq)^{-1}\big)^*\wh \clq^*J\clq  (\mfa_2(z)\clq)^{-1}.
\end{align} 
Setting 
\begin{align} &
 \label{Ve18!} 
\wh \clq= \clq, \quad \wh \vp=\vp,
\end{align}  
and taking into account that $\om_rJ\om_r^*=0$, we see that \eqref{Ve19} implies \eqref{Ve18}.

Formula \eqref{Ve18} shows that $\tau$ in \eqref{H10} is piecewise constant with
a finite number of jumps. Now, similar to the case \eqref{Ve14}, relations \eqref{H10} and \eqref{H11} yield
asymptotic formula
\begin{align} &
 \label{Ve20} 
\vp(z)=-\left(\sum_{k=0}^{2r-3}\frac{1}{z^{k+1}} H_k+\frac{1}{z^{2r-1}} \breve H_{2r-2}+\frac{1}{z^{2r}} \breve H_{2r-1}\right)+O\left(\frac{1}{z^{2r+1}}\right)
\end{align}  
for $z \to \infty$, and we set
\begin{align} &
 \label{Ve21} 
H_{2r-1}=\breve H_{2r-1}, \quad H_{2r}=(\I \om_{r}J \om_{r-1}^*)^{-1}+\clh_{21}(r+1)T(r)\clh_{12}(r+1).
\end{align} 

The inversion formula (of the form \eqref{Hf1}) for the extension $H(r+1)$ of $H(r)$ (defined via \eqref{Ve21})
is constructed using a triangular factorization of  $H(r+1)$ and it is easy to see that the invertibility of $H(r+1)$
follows from the invertibility of $H(r)$ and of $H_{2r}-\clh_{21}(r+1)T(r)\clh_{12}(r+1)$. Moreover, rewriting
\eqref{Hf1} for $T(r+1)=H(r+1)^{-1}$ we have
\begin{align}& \label{Ve22}
T(r+1)=\begin{bmatrix} T(r) & 0
\\ 0 & 0 \end{bmatrix}+
\begin{bmatrix} T(r)\clh_{12}
\\ -I_p
 \end{bmatrix}t_{r+1}
\begin{bmatrix} \clh_{12}^*T(r) & -I_p \end{bmatrix}>0,
\end{align}
and so $H(r+1)=T(r+1)^{-1}>0$. Thus, $H(r+1)$ generates a sequence of Verblunsky-type coefficients
$\om_0, \ldots, \om_{r-1}, \wt \om_r$ and it remains to show that
\begin{align} &
 \label{Ve23} 
\wt \om_r=\om_r,
\end{align}  
where $\om_r$ was given in the conditions of the proposition and was applied to the construction of $H(r+1)$ in 
\eqref{Ve20}, \eqref{Ve21}.

Since $\wt \om_r$ is a Verblunsky-type coefficient, it follows from Proposition \ref{PnProp} that
$-\I \om_{r-1} J \wt \om_r^*=t_{r+1}$. On the other hand, from \eqref{Hf0++} and \eqref{Ve21} we derive
$-\I \om_{r-1} J  \om_r^*=t_{r+1}$. Hence, the formula below is valid:
\begin{align} &
 \label{Ve23!} 
-\I \om_{r-1} J \wt \om_r^*=t_{r+1}=-\I \om_{r-1} J  \om_r^*.
\end{align}  

 Setting  $n=r+1$  and $\clq=\wt \clq\equiv  \wt \om_r^*$  in Theorem \ref{Hamburg},
 we construct via \eqref{H8} a function $\wt \vp(z)$. Similar to \eqref{Ve18} we have the equality
 \begin{align} &
 \label{Ve23+} 
\wt \vp(z)=\wt \vp(\ov{z})^*.
\end{align}  
 Similar to \eqref{Ve14} and \eqref{Ve20}, the function $\wt \vp$ has the
 asymptotic expansion 
 \begin{align}
 \label{Ve24} &
\wt \vp(z)=-\sum_{k=0}^{2r-1}\frac{1}{z^{k+1}} H_k+O\left(\frac{1}{z^{2r+1}}\right).
\end{align}  
According to \eqref{H7+}, \eqref{Hf4}, \eqref{Ve1} and \eqref{Ve3}, we may factorize $\mfa(r+1,z)$
from the representation \eqref{H8} of $\wt \vp$ in the following way:
\begin{align}
 \label{Ve25} &
\mfa(r+1,z)=\mfa(r,z)(I_{2p}-\I z \wt \om_r^* t_{r+1}^{-1} \wt \om_r J).
\end{align}  
Taking into account the property $\wt \om_r J \wt \om_r^*=0$ of the Verblunsky-type coefficient $\wt \om_r$
and the equality $\wt \clq\equiv  \wt \om_r^*$, we obtain from \eqref{Ve25} the relation
\begin{align}
 \label{Ve26} &
\mfa(r+1,z)\wt \clq=\mfa(r,z)\wt \clq.
\end{align}  
Therefore, we set $\wh Q= \wt Q$ and rewrite the representation of $\wt \vp$ of the form \eqref{H8}
in terms of $\mfa(r,z)$ instead of $\mfa(r+1,z)$:
\begin{align}& \label{Ve27}
\wt \vp(z) =\wh \vp(z)=\wh \vp(\ov{z})^*=\I \mfa_1(z)\wh \clq\big(\mfa_2(z)\wh \clq\big)^{-1}, \\
& \label{Ve28}
 \mfa_1(z):=\begin{bmatrix}  I_p & 0 \end{bmatrix}\mfa(r,z),
 \quad \mfa_2(z):=\begin{bmatrix}  0 & I_p \end{bmatrix}\mfa(r,z), \quad \wh \clq=\wt \clq= \wt \om_r^*.
\end{align}
We note that $\mfa(r,z)=\prod_{k=0}^{r-1}((I_{2p}-\I z Q_k J)$, and  so in view of  \eqref{Ve3}, \eqref{Ve13} and \eqref{Ve13+}
the leading terms of the matrix polynomials $\mfa_2(z)\clq$ and $\mfa_2\wh \clq$
have the form
\begin{align}
 \label{Ve29} &
\mfa_2(z)\clq=(-\I \om_r J \clq)z^r+O(z^{r-1}), \quad \mfa_2(z)\wh \clq=(-\I \om_r J \wh \clq)z^r+O(z^{r-1}).
\end{align}  
Using \eqref{Ve17}, \eqref{Ve23!} and the third equality in \eqref{Ve28}, we rewrite \eqref{Ve29} as:
\begin{align}
 \label{Ve30} &
\mfa_2(z)\clq=t_{r+1}z^r+O(z^{r-1}), \quad \mfa_2(z)\wh \clq=t_{r+1}z^r+O(z^{r-1}).
\end{align}  
According to \eqref{Ve19}, \eqref{Ve27} and \eqref{Ve30}, we have
\begin{align}
 \label{Ve31} &
\vp(z)-\wt \vp(z)= \I z^{-2r}\big(t_{r+1}^{-1}\wh \clq^* J \clq t_{r+1}^{-1}\big)+O(z^{-2r-1}).
\end{align}  
On the other hand, from \eqref{Ve20}, \eqref{Ve21} and \eqref{Ve24} we derive
\begin{align}&
 \label{Ve32} 
\vp(z)-\wt \vp(z)= z^{-2r+1}\big( H_{2r-2}-\breve H_{2r-2}\big)+O(z^{-2r-1}).
\end{align}
Compare \eqref{Ve31} and \eqref{Ve32} in order to see that
\begin{align}
 \label{Ve33} &
H_{2r-2}=\breve H_{2r-2}, \quad \wh \clq^* J \clq =0.
\end{align}
Formula \eqref{Ve17} and the last equalities in \eqref{Ve28} and \eqref{Ve33} show that
\begin{align}&
 \label{Ve34} 
 \wt \om_r J \om_r^* =0.
\end{align}
Relations \eqref{Ve16}, \eqref{Ve16+}, \eqref{Ve23!}, \eqref{Ve34} and $\wt \om_r J \wt \om_r^*=0$
imply the following equalities
\begin{align}&
 \nn 
\clk J\clk^*=\begin{bmatrix}0 & t_{r+1} \\ t_{r+1} & 0 \end{bmatrix}, \quad \clk \om_r^*=\clk  \wt \om_r^*=\begin{bmatrix}0 \\ t_{r+1} \end{bmatrix}
\quad {\mathrm{for}} \quad \clk :=\begin{bmatrix}\om_r J \\ -\I \om_{r-1} J \end{bmatrix}.
\end{align}
In view of the first equality in the formula above, we have $\det \clk \not=0$, and so the second equality yields \eqref{Ve23}.
\end{proof}
The first equality in \eqref{Ve33} means that for $\clq=\om_r^*$ the inequality in \eqref{H12}, where $n$ is substituted for $r$, turns into equality.
Note  also that if $\clq$ satisfies \eqref{H6} (i.e., $\clq$ is nonsingular with property-$J$), then the matrix $\clq u$, where $u$ is a $p \times p$ matrix and $\det u\not=0$,
is nonsingular with property-$J$ as well. Moreover, $\vp$ in \eqref{H8} does not change if we substitute $\clq$ by $\clq u$, that is, $\vp_{\clq}=\vp_{\clq u}$.
Hence, we may talk about  $\clq=\om_r^*u$ and have the following corollary.
\begin{Cy} When $\clq$ in Theorem \ref{Hamburg} has the properties
\begin{align}&
 \label{Ve35} 
 \clq^* J\clq=0, \quad \det(\om_{n-1}J \clq)\not=0,
\end{align}
this $\clq$ is nonsingular with property-$J$ and  the equality
\begin{align}& \label{H12m}
H_{2n-2} = \int_{-\infty}^{\infty}t^{2n-2}d\tau(t) \end{align}
is valid.
\end{Cy} 
Propositions \ref{PnUniqC} and \ref{PnOtO} imply the Verblunsky-type theorem below.
\begin{Tm} \label{TmM2} There is a one to one correspondence between Hankel matrices $H(n)>0$ $(n\in \BN)$
and their Verblunsky-type coefficients. 

The blocks of $H(n)$ are recovered from the
Verblunsky-type coefficients of $H(n)$ successively.  Namely, $H_0$ is recovered from \eqref{H0},
the block $H_{2r-1}$ $(r>0)$ is uniquely recovered from the asymptotic expansions \eqref{Ve14}
of the matrix functions $\vp(z)=\I \mfa_1(r+1,z)\clq(z)\big(\mfa_2(r+1,z)\clq(z)\big)^{-1}$,
and the matrix $H_{2r}$ is recovered $($after we have $H_{2r-1})$ from \eqref{Ve15}.

For each fixed $n\in \BN$, the set of the sequences $\om_k$ $(0 \leq k<n)$ of Verblunsky-type coefficients
of Hankel matrices $H(n)>0$ coincides  with the set of sequences of matrices $\om_k$ 
satisfying  \eqref{Ve16} and \eqref{Ve16+}.
\end{Tm}

\paragraph{3.}  Now, we will consider interconnections between canonical systems
generated by Hankel matrices $H(n)$ and Verblunsky-type coefficients of these matrices
in greater detail.
\begin{Pn}\label{PnInterCV1}
For each fixed $n\in \BN$, the set of canonical systems \eqref{Ve1}--\eqref{Ve2+} generated by the Hankel matrices $H(n)>0$
coincides with the set of canonical systems
\begin{align}& \label{Ve36}
y_{k+1}=w_{k+1}(\la)y_{k},  \quad w_{k+1}( \la)=I_{2p}+\frac{\I}{\la}J Q_{k}, \quad Q_k=\g_k^*\g_k,
\end{align}
where $\g_k\in \BC^{p \times 2p}$ $(0\leq k<n)$ and
\begin{align} & \label{Ve37} 
 \g_0 \begin{bmatrix} I_p \\ 0\end{bmatrix}=0, \quad \g_kJ\g_k^*=0, \quad \det \g_{k-1}J\g_k^* \not=0 \quad (0<k<n).
\end{align}
\end{Pn}
\begin{proof}. According to \eqref{Ve3}, \eqref{Ve13} and \eqref{Ve13+}, each canonical system \eqref{Ve1}--\eqref{Ve2+} generated by some Hankel matrix $H(n)>0$
has the form \eqref{Ve36}, \eqref{Ve37}, where $\g_k=t_{k+1}^{-1/2} \om_k$ $(0\leq k<n)$. 

On the other hand, for each system \eqref{Ve36}, \eqref{Ve37} we may
choose successively (for the increasing values of $k$) such pairs of unitary matrices $u_k$ and $p\times 2p$ matrices $\om_k$ that
\begin{align}& \label{Ve38}
u_0\g_0 \begin{bmatrix}0  \\ I_p \end{bmatrix}>0, \quad \om_0:=\left(u_0\g_0 \begin{bmatrix}0  \\ I_p \end{bmatrix}\right)u_0\g_0;
\\ & \label{Ve39}
\I u_k \g_k J\om_{k-1}^*>0, \quad \om_k:=\left( \I u_k \g_k J\om_{k-1}^*\right)u_k\g_k  \quad (0<k<n).
\end{align}
Clearly, the existence of the corresponding matrices $u_k$  follows from the third relation in \eqref{Ve37} and from the existence
of polar decompositions of square complex matrices. Taking into account \eqref{Ve38} and \eqref{Ve39}, we obtain
\begin{align} & \label{Ve40} 
\om_0\begin{bmatrix}0  \\ I_p \end{bmatrix}=\left(u_0\g_0 \begin{bmatrix}0  \\ I_p \end{bmatrix}\right)^2>0, \quad \I \om_kJ\om_{k-1}^*=
\left( \I u_k \g_k J\om_{k-1}^*\right)^2>0
\end{align}
for $k>0$. Formulas \eqref{Ve37}--\eqref{Ve40} imply that 
the conditions \eqref{Ve16} and \eqref{Ve16+} hold, that is, by virtue of Proposition \ref{PnOtO} the matrices $\{\om_k\}$
are Verblunsky-type coefficients. Moreover, in view of \eqref{Ve38}--\eqref{Ve40}, the Hamiltonian $Q_k=\g_k^*\g_k$
coincides with the matrices $Q_k$ in \eqref{Ve8}. Hence, Corollary \ref{CyQ} yields that the system 
\eqref{Ve36}, \eqref{Ve37} is generated by the Hankel matrix $H(n)>0$ determined via the constructed above
Verblunsky-type coefficients $\{\om_k\}$.
\end{proof}
In the proof of Proposition \ref{PnInterCV1},  we have shown that each canonical system \eqref{Ve36}, \eqref{Ve37}
is determined via a sequence of Verblunsky-type coefficients using relations \eqref{Ve8}.
Let us prove the uniqueness of the sequence of Verblunsky-type coefficients, which determines any fixed system \eqref{Ve36}, \eqref{Ve37}.
Indeed, let Verblunsky-type coefficients $\{\om_k\}$ determine system \eqref{Ve36}, \eqref{Ve37}. 
Then, comparing equalities in \eqref{Ve8} with $Q_k=\g_k^*\g_k$, for some matrices $q_k$ we have
\begin{align} & \label{Ve41} 
\om_k=q_k \g_k, \quad \det q_k\not=0.
\end{align}
In particular, we have 
\begin{align} & \label{Ve42} 
\om_0=\begin{bmatrix}0 & \vk\end{bmatrix} \quad (\vk>0), \quad \g_0=\begin{bmatrix}0 & \chi\end{bmatrix} \quad (\det \chi \not=0).
\end{align}
Hence, \eqref{Ve8} and \eqref{Ve41} for $k=0$ imply $q_0\chi=\chi^*\chi$, and $q_0$ is uniquely defined.
Assume that the matrices $q_0,\ldots, q_{r-1}$ are uniquely defined
and consider $Q_{r}=\g_{r}^*\g_{r}$. The formulas \eqref{Ve8}  and \eqref{Ve41} for $k=r$
yield now
$$q_{r}^*(\I q_{r}\g_{r}J\om_{r-1}^*)^{-1}q(r)=I_p, \quad {\mathrm{e.g.}}, \quad q(r)^*=\I \g_{r}J\om_{r-1}^*.$$
Thus, $q(r)$ and $\om_r$ are uniquely defined. We proved by induction the following corollary.
\begin{Cy}\label{CyL}
There is a one to one correspondence between Verblunsky-type coefficients and canonical
systems \eqref{Ve36}, \eqref{Ve37}. This correspondence is given by the formula \eqref{Ve8}.
\end{Cy}
\begin{Rk}\label{RkL} It follows from Corollary \ref{CyL} that there is a one to one correspondence
between infinite sequences $\{\om_k\}$ $(0\leq k <\infty)$ satisfying relations \eqref{Ve16} and \eqref{Ve16+}, where $n$ is substituted
by $\infty$ $($i.e.,  infinite sequences of Verblunsky-type coefficients$)$, and  canonical systems \eqref{Ve36}, \eqref{Ve37} on the semiaxis $k\geq 0$
$($i.e., systems \eqref{Ve36}, \eqref{Ve37} with $n$  substituted in \eqref{Ve37}
by $\infty)$. It follows also from Theorem \ref{TmM2} that there is a one to one correspondence between infinite sequences 
of the blocks $\{H_k\}$, such that all the matrices $H(n)=\{H_{i+j-2}\}_{i,j=1}^n$ $(n\geq 0)$ are positive-definite,  and  infinite sequences
of Verblunsky-type coefficients.
\end{Rk}
%%%%%%%%%%%%%%%%%%%
\paragraph{4.} 
Consider system \eqref{Ve36}, \eqref{Ve37} on the semiaxis $k\geq 0$.
In order to establish connections with the spectral theory,
we will need the definition of the spectral function of  the discrete canonical system from  \cite[Section 8]{SaL2}.
For this purpose, we introduce the spaces $\ell^2_Q$ and $L^2(\tau)$, where $Q(k)=Q_k\geq 0$ and $\tau$
is a non-decreasing $p \times p$ matrix function on $\BR$.
The scalar product in $\ell^2_Q$  is given by the formula $(h, \wt h)_{Q}=\sum_{k=0}^{\infty}\wt h(k)^*Q_kh(k)$
and the scalar product in $L^2(\tau)$ is given by the formula $(f, \wt f)_{\tau}=\int_{-\infty}^{\infty}\wt f(t)^*d\tau(t)f(t)$.
Denote the normalized fundamental solution of the system \eqref{Ve36}, \eqref{Ve37} ($k\geq 0$) by $Y$, that is, let $Y$
satisfy relations
\begin{align}& \label{Ve43}
Y(k+1,\la)-Y(k,\la)=\frac{\I}{\la}J Q_{k}Y(k,\la) \quad (k\geq 0),  \quad Y(0)=I_p.
\end{align}
Finally, introduce the operator $V$ from $\ell^2_Q$ into $L^2(\tau)$ by its action on the
functions with finite support:
\begin{align}& \label{Ve44}
\big(Vh\big)(\la)=\begin{bmatrix}0 & I_p\end{bmatrix}\sum_{k=0}^{\infty}Y\Big(k, \frac{1}{\la}\Big)^*Q_kh(k).
\end{align}
\begin{Dn} A  non-decreasing $p \times p$ matrix function $\tau$ on $\BR$ is called a spectral function
of the system \eqref{Ve36}, \eqref{Ve37} on the semiaxis $k\geq 0$ if $V$ given in \eqref{Ve44} is an isometric mapping into 
$L^2(\tau)$.
\end{Dn}
Although we restricted the definition above to the system, which we consider here, it is clear
that the same definition  works for all canonical systems.

Slightly rephrasing Corollaries 8.2.1 and 8.2.3 from \cite{SaL2} we obtain the following theorem.
\begin{Tm}\label{TmF} A  non-decreasing $p \times p$ matrix function $\tau$ on $\BR$ is  a spectral function
of some system \eqref{Ve36}, \eqref{Ve37} on the semiaxis $k\geq 0$ if and only if all the blocks given by
\begin{align}& \label{Ve45}
H_k=\int_{-\infty}^{\infty}t^{k}d\tau(t) \quad (k \geq 0)
\end{align}
are well-defined and all the matrices $H(n)=\{H_{i+j-2}\}_{i,j=1}^n$ $(n\geq 1)$ are positive-definite.
Moreover, if  $\tau$  is, indeed,  a spectral function
of some system \eqref{Ve36}, \eqref{Ve37}, it is a spectral function of only one system
\eqref{Ve36}, \eqref{Ve37} $(k \geq 0)$.
\end{Tm}
%%%%%%%%%%%%%%%%
\bigskip

\noindent{\bf Acknowledgments.}
 {This research   was supported by the
Austrian Science Fund (FWF) under Grant  No. P29177.}

\appendix
%%%%%%
%%%%%%

\section{Appendix}
\begin{proof} of Theorem \ref{Hamburg}.
When looking at this proof, one has to take into account that $H$ is a particular case of operators $S$   in \cite{SaL2},
and that the notations $\b$ and $\a$ are used  in \cite{SaL2} instead of $\mu$ and $\nu$ in the
Herglotz representation \eqref{H14}.
Since 
\begin{align}& \label{H13}
H>0; \quad \Phi_2g\not=0 \quad {\mathrm{for \, \, all}} \quad g\not=0 \quad (g\in \BC^p); \quad \s(A)=0
 \end{align}
(where $\s(A)$ means the spectrum of $A$), the conditions of  \cite[Theorem 1.4.2]{SaL2} are fulfilled.
Thus, the set of functions (matrix functions) $\vp(z)$ constructed in \eqref{H8} coincides with the set of holomorphic (in $\BC_+$) solutions
of  the fundamental Potapov's inequality \cite[(1.2.1)]{SaL2}. 

In particular, the inequality 
$$(\vp(z)-\vp(z)^*)/(z- \ov{z}) \geq 0$$
is immediate from \cite[(1.2.1)]{SaL2}, that is, the functions  $\vp(z)$ belong to the Herglotz class. Hence, the functions $\vp$
admit Herglotz representation
\begin{align}& \label{H14}
\vp(z)=\mu z +\nu+\int_{-\infty}^{\infty}\left(\frac{1}{t-z}-\frac{t}{1+t^2}\right)d\tau(t) \quad (\mu \geq 0, \quad \nu=\nu^*).
\end{align}
Moreover, according to  \cite[Theorem 1.3.1]{SaL2} the nondecreasing $p \times p$ matrix functions $\tau$ in \eqref{H14} belong to $\cle(n)$.
Hence, the right-hand sides of \eqref{H11} and \eqref{H12} are well-defined. The inequality \eqref{H12} itself follows from the inequality 
\begin{align}& \label{H15}
H \geq  H_{\tau}:=\int_{-\infty}^{\infty}(I_{np}-tA)^{-1}\Phi_2 d\tau(t) \Phi_2^*(I_{np}-tA^*)^{-1}
\end{align}
given on \cite[p. 7] {SaL2} and from the easy relation
\begin{align}& \label{H16}
(I_{np}-zA)^{-1}\Phi_2 =\begin{bmatrix} I_p \\ z I_p \\ \cdots \\ z^{n-1} I_p \end{bmatrix}.
\end{align}

Recall that $\vp$ satisfies the  Potapov's inequality \cite[(1.2.1)]{SaL2}. Therefore, we may use \cite[Lemma 1.3.1]{SaL2},
and, more precisely, the equality
\begin{align}& \label{H17}
\Phi_2\mu \Phi_2^*=A (H-H_{\tau}) A^*.
\end{align}
In view of the definitions \eqref{H3} and \eqref{H4} of $A$ and $\Phi_2$, formula \eqref{H17} yields
\begin{align}& \label{H18}
\mu =0.
\end{align}
Taking into account the definition of $H_{\tau}$ in \eqref{H15} and relations \eqref{H16} and \eqref{H18}, we see that the same
formula \eqref{H17} implies \eqref{H11}.

Finally, from  \cite[Proposition 1.3.1]{SaL2} we obtain
\begin{align}& \label{H19}
\nu=\int_{-\infty}^{\infty}\frac{t}{1+t^2}d\tau(t).
\end{align}
Substituting \eqref{H18} and \eqref{H19} into \eqref{H14}, we derive \eqref{H10}.

We have shown that for all $\vp$ constructed in \eqref{H8} relations \eqref{H10}--\eqref{H12} hold.
It remains to prove that all nondecreasing functions $\tau$ satisfying \eqref{H11} and \eqref{H12}
are described via the set of  $\vp$ constructed in \eqref{H8}.

Indeed, assume that \eqref{H11} and \eqref{H12} hold for some nondecreasing $\tau(t)$, and define $\mu$ and $\nu$
via \eqref{H18} and \eqref{H19}. Then, \cite[Lemma 1.1.2]{SaL2} together with \eqref{H11} and \eqref{H12}  shows that \cite[(1.1.11)]{SaL2} 
is valid. Hence, by virtue of \cite[Corollary 1.2.1]{SaL2}, the function $\vp(z)$ determined by our $\tau$ via \eqref{H10} satisfies Potapov's inequality
\cite[(1.2.1)]{SaL2}. As discussed at the beginning of the proof, this means that $\vp$ is one of the functions constructed in \eqref{H8},
that is, $\tau$ belongs to the set given by \eqref{H8}--\eqref{H10}.
\end{proof}

%%%%%%%%%%%%%%%%%%%%%%%%%%%%%%%%%%%%%%%%%%%%%%
%%%%%%%%%%%%%%%%%%%%%%%%%%%%%%%%%%%%%%%%%%%%%%%
%\newpage
%%%%%%%%%%%%%%%%%%%%%%%%%%%%%%%%%%%%%%%%%%%
%%%%%%%%%%%%%%%%%%%%%%%%%%%%%%%%%%%%%%%%%%%%

\begin{flushright}

A.L. Sakhnovich,\\
Fakult\"at f\"ur Mathematik, Universit\"at Wien, \\
Oskar-Morgenstern-Platz 1, A-1090 Vienna, Austria\\
E-mail: oleksandr.sakhnovych@univie.ac.at
\end{flushright}

%%%%%%%%%%%%%%%%%%%%%%%%

\end{document}